\documentclass[12pt]{amsart}

 \newtheorem{thm}{Theorem}[section]
 
 \newtheorem{lem}[thm]{Lemma}
 \newtheorem{prop}[thm]{Proposition}
 \theoremstyle{definition}
 \newtheorem{defn}[thm]{Definition}
 \theoremstyle{remark}
 \newtheorem{rem}[thm]{Remark}
 \theoremstyle{definition}
 \newtheorem{ex}[thm]{Example}

 \newcommand{\PP}{\mathbb{P}}

\usepackage{xypic} \usepackage{amsthm}

\begin{document}

\title{Bipolynomial Hilbert functions}

\author[E. Carlini]{Enrico Carlini}
\address[E. Carlini]{Dipartimento di Matematica, Politecnico di Torino, Torino, Italia}
\email{enrico.carlini@polito.it}

\author[M.V.Catalisano]{Maria Virginia Catalisano}
\address[M.V.Catalisano]{DIPTEM - Dipartimento di Ingegneria della Produzione, Termoenergetica e Modelli
Matematici, Universit\`{a} di Genova, Piazzale Kennedy, pad. D
16129 Genoa, Italy.} \email{catalisano@diptem.unige.it}

\author[A.V. Geramita]{Anthony V. Geramita}
\address[A.V. Geramita]{Department of Mathematics and Statistics, Queen's University, Kingston, Ontario, Canada, K7L 3N6 and Dipartimento di Matematica, Universit\`{a} di Genova, Genova, Italia}
\email{Anthony.Geramita@gmail.com \\ geramita@dima.unige.it  }


\begin{abstract}
Let $X\subset\PP^n$ be a closed subscheme and let $HF(X,\cdot)$
and $hp(X,\cdot)$ denote, respectively, the Hilbert function and
the Hilbert polynomial of $X$. We say that $X$ has {\it
bipolynomial Hilbert function} if
$HF(X,d)=\min\left\{hp(\PP^n,d),hp(X,d)\right\}$ for every
$d\in\mathbb{N}$. We show that if $X$ consists of a plane and
generic lines, then $X$ has bipolynomial Hilbert function. We also
conjecture that generic configurations of non-intersecting linear
spaces have bipolynomial Hilbert function.
\end{abstract}

\maketitle

\section{Introduction}

The Hilbert function of a scheme $X\subset\PP^n$ encodes a great
deal of interesting information about the geometry of $X$ and so
the study of  $HF(X,\cdot)$ has generated an enormous amount of
research. One of the most crucial and basic facts about the
Hilbert function of a scheme is that the function is eventually
polynomial. More precisely
\[HF(X,d)=hp(X,d), \mbox{ for } d\gg 0.\]

In general, knowledge of the Hilbert polynomial does not determine
the Hilbert function. But, there are some interesting situations
when this is the case. E.g. if $X$ is a generic set of $s$ points
in $\PP^n$, it is well known, and not hard to prove, that
\[HF(X,d)=\min\left\{hp(\PP^n,d)={n+d\choose d},hp(X,d)=s\right\},\]
for all $d\in\mathbb{N}$. A much harder result is due to
Hartshorne and Hirschowitz. In \cite{HartshorneHirschowitz} the
authors considered schemes $X\subset\PP^n$ consisting of $s$
generic lines and they proved that
\[HF(X,d)=\min\left\{hp(\PP^n,d)={n+d\choose d},hp(X,d)=s(d+1)\right\},\]
for all $d\in\mathbb{N}$.

Inspired by these results about points and  lines, we restrict our
attention to that special family of schemes known as {\em
configurations of linear spaces}. We recall that a configuration
of linear spaces $\Lambda\subset \PP^n$ is nothing more than a
finite collection of linear subspaces of $\PP^n$; see
\cite{CarCatGer1,CaCat09} and \cite{DerksenSidman} for more on
these schemes and their connection with {\em subspace
arrangements}. We further say that a configuration of linear
spaces is {\em generic} when its components are generically
chosen.

The Hilbert polynomial of a generic configuration of linear spaces
is known, thanks to a result of Derksen, see \cite{Derksen}. Thus,
in light of the results on the Hilbert function of generic points
and generic lines, we propose the following
\begin{quote} {\bf Conjecture:} {\em if $\Lambda\subset \PP^n$ is a
generic configuration of linear spaces with non-intersecting
components, then
\[HF(X,d)=\min\left\{hp(\PP^n,d),hp(X,d)\right\},\]
for all $d\in\mathbb{N}$.}
\end{quote}
We will call a Hilbert function defined as above {\em
bipolynomial}. Hence, the conjecture states that generic
configurations of linear spaces with non-intersecting components
have bipolynomial Hilbert function.

As we mentioned above, this conjecture is true when
$\dim\Lambda=0$ (generic points) and when $\dim\Lambda=1$. The
conjecture holds in the dimension one case because of the result
about generic lines in \cite{HartshorneHirschowitz} and because
we know how adding generic points to a scheme
changes its Hilbert function, see \cite{GeramitaMarosciaRoberts}.

In this paper we produce new evidence supporting our conjecture.
Namely, we show that the union of one plane and $s$ generic lines
 has bipolynomial Hilbert function.

The paper is structured as follows: in Section \ref{basicsection}
we introduce some basic notation and results we will use;
Sections \ref{p4section} and \ref{lemmas} contain the base cases for our inductive approach;
Section \ref{generalsection} contains our main result, Theorem
\ref{TeoremaInPn}.  These sections are followed by a section on Applications and another in which we propose a possibility for the Hilbert function of any generic configuration of linear spaces, even one in which there are forced intersections.

The first two authors thank Queen's University for its hospitality during part of the
preparation of this paper.  All the authors enjoyed support from
 NSERC (Canada)  and   GNSAGA of INDAM (Italy). The first author was, furthermore,
partially supported by a ``Giovani ricercatori, bando 2008" grant
of the Politecnico di Torino.

\section{Basic facts and notation}\label{basicsection}

We will always work over an algebraically closed field $k$ of
characteristic zero. Let $R=k[x_0,...,x_n]$ be the coordinate ring
of $\PP^n$, and denote by  $I_X$ the ideal of a scheme $X \subset
\PP^n$. The Hilbert function of $X$ is then  $HF(X,d)=\dim
(R/I_X)_d$.

\begin{defn} \ Let $X$ be a subscheme of $\PP^n$.  We say that $X$ has a {\it bipolynomial Hilbert function} if
\[HF(X,d)=\min\left\{hp(\PP^n,d),hp(X,d)\right\},\]
for all $d\in\mathbb{N}$.
\end{defn}

It will often be convenient to use ideal notation rather then Hilbert function notation, i.e. we will often describe $\dim I_X$ rather than $HF(X,d)$. It is clearly trivial to pass
from one piece of information to the other.

The following lemma gives a criterion for adding to a scheme, $X
\subseteq \Bbb P^ n$, a set of reduced points lying on a linear
space ${\Pi} \subseteq \Bbb P ^n$ and imposing independent
conditions to forms of a given degree in the ideal of $X$.

 \begin{lem} \label{AggiungerePuntiSuSpazioLineare}
Let $d \in \Bbb N$. Let $X  \subseteq \Bbb P^n$ be  a scheme, and
let $P_1,\dots,P_s$ be generic distinct points on a  linear space
$\Pi \subseteq \Bbb P^n$.

If $\dim (I_{X })_d =s$ and $\dim (I_{X +\Pi})_d =0$, then
$
\dim (I_{X+P_1+\cdots+P_s})_d = 0.$
\par

\end{lem}

\begin{proof}

By induction on $s$. Obvious for  $s=1$. Let $s>1$ and let $X' =
X+P_s$. Obviously $\dim (I_{X ' + \Pi})_d =0$. Since $\dim (I_{X
+\Pi})_d=0$ and $P_s$ is a generic point in $\Pi$, then $\dim
(I_{X' })_d =s-1$. Hence, by the inductive hypothesis, we get
$\dim (I_{X'+P_1+\cdots+P_{s-1}})_d =\dim (I_{X+P_1+\cdots+P_s})_d
= 0.$

\end{proof}

Since we will make use of Castelnuovo's inequality several times in the
next sections, we recall it here in a form more suited to our use (for
notation and proof we refer to \cite{AH95}, Section 2).

\begin {defn}\label{ResiduoTraccia}
If $X, Y$ are closed subschemes of $\PP^n$, we denote by $Res_Y X$
the scheme defined by the ideal $(I_X:I_Y)$ and we call it the
{\it residual scheme} of $X$ with respect to $Y$, while the scheme
$Tr_Y X \subset Y$ is the schematic intersection $X\cap Y$, called
the {\it trace} of $X$ on $Y$.

\end {defn}

 \begin{lem} \label{Castelnuovo}{\bf (Castelnuovo's inequality):}
Let $d,\delta \in \mathbb N$, $d \geq \delta$, let ${Y} \subseteq \PP ^n$ be a smooth hypersurface of degree $\delta$,
and let $X \subseteq \PP
^n$ be a  scheme. Then
$$
\dim (I_{X, \PP^n})_d  \leq  \dim (I_{ Res_Y X, \PP^n})_{d-\delta}+
\dim (I_{Tr _{Y} X, Y})_d.
$$
\qed
\end{lem}

Even though we will only use the following lemma in the cases
$m=2$, $m=3$ (see the notation  in the lemma), it seemed
appropriate to give the more general argument since such easily
understood (and non trivial) degenerations occur infrequently.
 \begin{lem} \label{sundial}
 Let $X_1  \subset \Bbb P^n $ be the disconnected  subscheme consisting of a  line $L_1$ and a linear space $\Pi \simeq \mathbb P^m$ (so the linear span of $X_1$ is $<X_1> \simeq \Bbb P^{m+2} $). Then there exists a flat family  of subschemes $$X_{\lambda}\subset <X_1> \ \ \ \ \ (\lambda \in k )$$
whose special fibre $X_0$ is the union of

\begin{itemize}
\item
the linear space $\Pi $,

\item a line $L$ which intersects  $\Pi$ in a point $P$,

\item   the scheme $2P|_ { <X_1>}$, that is, the schematic intersection of the double point
$2P$ of $\mathbb P^n$ and $<X_1>$.

\end{itemize}
Moreover, if $H \simeq  \Bbb P^{m+1}$ is the linear span of $L$
and  $\Pi$, then $Res_H(X_0)$ is given by the (simple) point $P$.

 \end{lem}

\begin{proof} We may assume that the ideal of the line $L_1$ is
 $$(x_1, \dots, x_{m}, x_{m+1} - x_{0} , x_{m+3},  \dots, x_n)$$ and the ideal of $\Pi$ is
$( x_{m+1}  , \dots, x_n)$, so the ideal of $X_1$ is
$$I_{X_1}=(x_1, \dots, x_{m}, x_{m+1} - x_{0} , x_{m+3},  \dots, x_n) \cap ( x_{m+1}  , \dots, x_n).
$$
Consider the flat family $\{ X_{\lambda}\}_{\lambda \in k}$, where for any fixed $\lambda \in k$, $X_\lambda$ is the union of  $\Pi$ and the line
$$
x_1= \dots= x_{m}= x_{m+1} - \lambda x_{0} = x_{m+3}= \dots= x_n=0.$$

The ideal of  $X_{\lambda}$ is
$$I_{X_{\lambda}}=(x_1, \dots, x_{m}, x_{m+1} - \lambda x_{0} , x_{m+3},  \dots, x_n)  \cap ( x_{m+1}  , \dots, x_n)
$$
$$=( x_1, \dots, x_{m}, x_{m+1}-\lambda x_0) \cap (x_{m+1},x_{m+2}) + (x_{m+3},\dots,x_{n})$$
$$=( x_1, \dots, x_{m}, x_{m+1}-\lambda x_0) \cdot (x_{m+1},x_{m+2}) + (x_{m+3},\dots,x_{n})$$
$$=(  x_1 x_{m+1}, \dots, x_{m} x_{m+1}, (x_{m+1}-\lambda x_0) x_{m+1})
$$
$$+
( x_1x_{m+2}, \dots, x_{m}x_{m+2}, (x_{m+1}-\lambda x_0)x_{m+2})$$
$$+ (x_{m+3},\dots,x_{n}),$$
which for $\lambda =0$ gives:
$$I_{X_0}=(  x_1 x_{m+1}, \dots, x_{m} x_{m+1}, x_{m+1}^2)+
( x_1x_{m+2}, \dots, x_{m}x_{m+2}, x_{m+1}x_{m+2})$$
$$+ (x_{m+3},\dots,x_{n}),$$
$$
=( x_1, \dots, x_{m+1}) \cdot (x_{m+1},x_{m+2}) + (x_{m+3},\dots,x_{n}).$$

Let $(x_{m+3},\dots,x_{n}) = J$.
We will prove that
\begin{equation} \label{idealediX0}
I_{X_0}
=( x_1, \dots, x_{m+1}) \cdot (x_{m+1},x_{m+2}) +J
\end{equation}
$$
=\left [( x_1, \dots, x_{m+1})+J\right ] \cap \left [(x_{m+1},x_{m+2})+J\right ] \cap
\left [ ( x_1, \dots,x_{m+2})^2 +J\right ]. $$

We use  Dedekind's Modular Law several times in what follows (see \cite[page 6]{AtMac}).
We start by considering the intersection of the first two ideals, i.e.,
$$
\left [( x_1, \dots, x_{m+1})+J\right ] \cap \left [(x_{m+1},x_{m+2})+J\right ] $$
$$=
\left [( x_1, \dots, x_{m+1})+J\right ] \cap \left [((x_{m+1})+J)+(x_{m+2})\right ]$$
$$=
 ((x_{m+1})+J)+
\left\{ \left [( x_1, \dots, x_{m+1})+J \right ] \cap (x_{m+2})\right \}$$
$$=
 ((x_{m+1},x_1x_{m+2}, \dots, x_mx_{m+2})+ J).$$

\medskip
 It remains to intersect this last ideal with the third ideal above, i.e.,
$$((x_{m+1},x_1x_{m+2}, \dots, x_mx_{m+2})+ J) \cap \left [ ( x_1, \dots,x_{m+2})^2 +J \right ]$$
$$=\left[(x_{m+1}) +((x_1x_{m+2}, \dots, x_mx_{m+2})+ J)\right] \cap \left [ ( x_1, \dots,x_{m+2})^2 +J \right ]$$
$$=\left[(x_{m+1})  \cap  (( x_1, \dots,x_{m+2})^2 +J )\right ] + ((x_1x_{m+2}, \dots, x_mx_{m+2})+ J)$$
$$=\left\{(x_{m+1})  \cap   \left[ (x_{m+1})\cdot ( x_1, \dots,x_{m+2}) + ( x_1, \dots,x_{m},x_{m+2})^2 +J
\right ] \right \}$$
$$+ ((x_1x_{m+2}, \dots, x_mx_{m+2})+J)$$
$$= \left[ (x_{m+1})\cdot ( x_1, \dots,x_{m+2}) \right ]+ \left [ (x_{m+1})  \cap \left (  ( x_1, \dots,x_{m},x_{m+2})^2 +J
\right ) \right ]
$$
$$+ ((x_1x_{m+2}, \dots, x_mx_{m+2})+ J).$$

Clearly the middle ideal is contained in the sum of the other two, and so the last ideal is equal to
$$\left[ (x_{m+1})\cdot ( x_1, \dots,x_{m+2}) \right ]+ ((x_1x_{m+2}, \dots, x_mx_{m+2})+ J)$$
$$=( x_1, \dots, x_{m+1}) \cdot (x_{m+1},x_{m+2}) +J.$$
So we have proved that $I_{X_0} $ is
$$\left [( x_1, \dots, x_{m+1})+J\right ] \cap \left [(x_{m+1},x_{m+2})+J\right ] \cap
\left [ ( x_1, \dots,x_{m+2})^2 +J\right ]. $$ Since $J $ is the
ideal of $<X_1>$, the first ideal in this intersection defines a
line $L$ in $<X_1>$ which meets the linear space $\Pi$ (defined by
the second ideal in this intersection) in the point
$P=[1:0:\dots:0]\in \mathbb P^n$, which is the support of the third
ideal in this intersection.  The  third ideal, in fact, describes
the scheme $2P|_ {<X_1>}$which is  the double point $2P$ of
$\mathbb P^n$  restricted to  the span of $X_1$.

The ideal of $H$ is  $(x_{m+1})+J $, hence from (\ref {idealediX0}) we have that the ideal of $Res_H(X_0)$ is
$$I_{X_0} : I_{H} = \left[ ( x_1, \dots, x_{m+1}) \cdot (x_{m+1},x_{m+2}) + J \right]
:((x_{m+1})+J) $$
$$=(x_1,\dots ,x_n)= I_P.
$$
\end{proof}

\begin {defn}\label{conica degenere}
We say that $C$ is a {\it degenerate conic} if  $C$ is the union
of two intersecting lines $L_{1}, L_{2}.$  In this case we write
$C=L_1+L_2$.
\end {defn}

\begin {defn}\label{defsundial}
Let  $n\geq m+2$. Let   $\Pi \simeq \Bbb P^m \subset \mathbb P^n$
be a linear space of dimension $m$, let $P\in \Pi$ be a point and
let $L \not\subset \Pi $ be a generic line through $P$. Let $T
\simeq \Bbb P^{m+2}$ be a generic linear space containing the
scheme $L+\Pi$.  We call the scheme $L+\Pi+ 2P|_T$ an {\it
$(m+2)$-dimensional  sundial}. (See, for instance,   the scheme
$X_0$ of Lemma \ref{sundial}).

\medskip
Note that for  $m=1$, the scheme $L+\Pi$ is a degenerate conic and
 the  $3-$dimensional  sundial
$L+\Pi+ 2P|_T$ is a {\it degenerate conic with an embedded point}
(see \cite {HartshorneHirschowitz}).
\end{defn}

\begin{thm} [Hartshorne-Hirschowitz, \cite{HartshorneHirschowitz}] \label{HH}
Let $n, d \in \mathbb N$.
For $n\geq 3$, the ideal of the scheme $X\subset \Bbb P^n$   consisting of $s$ generic
 lines has the expected dimension, that is,
$$
\dim (I_X)_d = \max \left \{ {d+n \choose n} -s(d+1), 0 \right \},
$$
or equivalently
$$
H(X,d) = \min \left \{ hp(\PP^n,d)={d+n \choose n}, hp(X,d)=s(d+1)
\right \}.
$$
 \end{thm}
\qed
\medskip

Since a line imposes at most $d+1$ conditions to the forms of degree $d$, the first part of the following lemma   is clear. The second statement of the lemma is obvious.

\begin{lem} \label{BastaProvarePers=e,e*}
Let $n,d, s  \in \Bbb N$,  $n \geq4$. Let $\Pi \subset \Bbb P^n$ be a plane, and let $L_1, \dots , L_s \subset \Bbb P^n$ be $s$ generic lines. Let
$$ X_s= \Pi +  L_1+ \dots + L_s \subset \Bbb P^n .$$
\begin{itemize}
\item[(i)] If  $ \dim (I_{X_s})_d = {d+n \choose n} -  {d+2 \choose 2}  -s(d+1)$, then
 $ \dim (I_{X_{s'}})_d = {d+n \choose n} -  {d+2 \choose 2}  -s'(d+1)$  for any $s'<s$.
 \par
\item[(ii)]  If  $ \dim (I_{X_s})_d = 0$, then
 $ \dim (I_{X_{s'}})_d =0$  for any $s'>s$.
 \par
\end{itemize}
\end{lem}
\qed

\section{The base for our induction}\label{p4section}

In this section we prove our main Theorem (see \ref{TeoremaInPn}) in $\PP^4$.

 \begin{thm} \label{TeoremaInP4} Let $d\in\mathbb{N}$ and $\Pi \subset \Bbb P^4$ be a plane,
 and let $L_1, \dots , L_s \subset \Bbb P^4$ be $s$ generic lines. Set
  $$X= \Pi +  L_1+ \dots + L_s \subset \Bbb P^4.$$
Then
$$
\dim (I_{X})_d = \max \left \{ {d+4 \choose 4} - {d+2
\choose 2} -s(d+1), 0 \right \},
$$
or equivalently $X$ has bipolynomial Hilbert function.

 \end{thm}

 \begin{proof} We proceed by induction on $d$. Since the theorem is obvious for $d=1$, let $d>1$. By  Lemma \ref{BastaProvarePers=e,e*}
 it suffices to prove the theorem for $s=e$ and $s=e^*$, where
 $$e= \left \lfloor   {{{d+4 \choose 4} -  {d+2 \choose 2} }\over {d+1} }\right \rfloor =
   \left \lfloor   {\frac{d (d+2)(d+7)}{24}   }\right \rfloor ; \ \ \ \
 e^*= \left \lceil   {{{d+4 \choose 4} -  {d+2 \choose 2} }\over {d+1} }\right \rceil .
 $$
   Let $$
 \bar e =  \left \lfloor {{{(d-1)+4 \choose 4} -  {(d-1)+2 \choose 2} }\over {(d-1)+1} } \right \rfloor =  \left \lfloor{{ (d-1)(d+1)(d+6)} \over 24} \right \rfloor .
 $$

We consider two cases.

 \par
 \medskip
 {\it Case 1}: $d$ odd.
 \par
For $s=e$,  we have to prove that $\dim (I_{X})_d =  {d+4
\choose 4} -  {d+2 \choose 2}  -e(d+1)$ (which is obviously positive).
Since $\dim (I_{X})_d \geq  {d+4 \choose 4} -  {d+2 \choose 2}  -e(d+1),$ we have
only to show that $\dim (I_{X})_d \leq  {d+4 \choose 4} -
{d+2 \choose 2}  -e(d+1).$

For $s=e^*$,  we have to prove that
$\dim (I_{X})_d =0.$

In order to prove these statements we construct a scheme $Y $
obtained from $X$ by
specializing the $s- \bar e$ lines $L_{\bar e +1}, \dots, L_{s}$
into a generic hyperplane $H \simeq \Bbb P^3$ (we can do this since $\bar e <s$).

 If we can prove that $\dim(I_{Y})_d =  \max \{ {d+4 \choose 4} -  {d+2 \choose 2}  -s(d+1); 0 \}$, that is,
 if we can show that  the plane and the $s$ lines give the expected number of conditions to the forms of degree $d$ of  $\mathbb P^4$,
then (by the semicontinuity of the Hilbert function)  we are done.
\par
Note that
 $$Res_H Y = L_{1}+ \dots + L_{\bar e} + \Pi \subset \Bbb P^4 ,
$$
and
$$Tr_H Y  =  P_{1}+ \dots + P_{\bar e }+ L_{\bar e +1}+ \dots+L_{s}+ L  \subset \Bbb P^3,$$
where $P_i = L_i \cap H$, ($1 \leq i \leq \bar e$), and $L$ is the line $\Pi \cap H$ .

\medskip
\medskip
Since $ d$ is odd, the number ${{ (d-1)(d+1)(d+6)} \over 24} $ is an integer, so
 $$
 \bar e ={{ (d-1)(d+1)(d+6)} \over 24}.
 $$
The inductive hypothesis applied to $Res_H Y $ in degree $d-1$ yields:
$$\dim (I_{Res_H Y })_{d-1} = {d+3 \choose 4} -  {d+1 \choose 2}  -\bar e (d) =0.
$$
By Theorem \ref{HH}, since the $P_i$ are generic points, we get
$$ \dim  (I_{ Tr_H Y  })_{d} =\max \left \{ {d+3 \choose 3} - \bar e - (s - \bar e +1) (d+1);0 \right \}$$
$$=
\max \left \{ {d+3 \choose 3} +\bar e d- (s  +1) (d+1) ;0\right \}$$

$$=\left\{
\begin{array}{cc}
{d+4 \choose 4} - {d+2 \choose 2} - e(d+1) & \ {\rm for} \  s=e \\
0 & \ \ \ {\rm for} \  s=e^* \\
 \end{array} \right. ,
$$
and the conclusion follows by Lemma \ref{Castelnuovo} with $\delta=1$.
\par

\medskip
 {\it Case 2}: $d$ even. \par
 In this case $e = e^* = {{d (d+2)(d+7)} \over 24} $, and so
we only have to prove that $\dim (I_{ X})_d =  0.$ Let
$$ x= {{d(d+2)}\over 8}  ,
$$
and note that $x$ is an integer,  $x <e$.

Let  $H \simeq \Bbb P^3$ be a generic hyperplane containing the
plane $\Pi$, and let  $ Y $ be the scheme obtained from $ X$ by
degenerating the  $x$ lines $L_{1}, \dots, L_{x}$ into  $H $. By
abuse of notation, we will again denote these lines by $L_{1},
\dots, L_{x}$. By Lemma \ref{sundial}, with $m=2$, we get
$$ Y = L_{1}+ \dots+ L_{x} + 2P_1+ \dots +2 P_x + \Pi + L_{x+1}+ \dots + L_{ e}  ,
$$
where $P_i= L_i \cap \Pi$ ($1 \leq i \leq x $) and the $2P_i$ are  double points in $\mathbb P^4$.
 If we can prove that $\dim (I_{ Y})_d =  0$ we are done.

By Lemma \ref{sundial}, with $m=2$, we get
 $$Res_H  Y = P_1+ \dots +P_x+ L_{x+1}+ \dots + L_{ e} \subset \Bbb P^4 ,
$$
where the $P_i$ are generic points in $\Pi$.

Also,
$$
Tr_H  Y  =   L_{1}+ \dots+ L_{x} + 2P_1|_ H+ \dots +2 P_x |_ H +
\Pi + Q_{x+1}+ \dots + Q_{ e}, $$ but, since $2P_i |_ H \subset
L_i + \Pi$, we get
$$
Tr_H  Y  =  L_{1}+ \dots+ L_{x} + \Pi  + Q_{x+1}+ \dots + Q_{ e}
\subset \Bbb P^3
$$
where $Q_i = L_i \cap H$, ($x+1 \leq i \leq  e$).

Since $\Pi$ is a fixed component of the zero locus for the forms
of $I_{ Y \cap H }$, we get that
$$ \dim  (I_{Tr_H  Y   })_{d} =\dim  (I_{Tr_H Y  - \Pi  })_{d-1}.
$$
Since the $Q_i$ are generic points, we can apply   Theorem \ref{HH} and get
\begin{equation} \label{traccia}
\dim  (I_{Tr_H Y  - \Pi  })_{d-1}={d-1+3 \choose 3} - xd - (e-x) =
0 .
\end{equation}
Now  we will prove that  $\dim (I_{Res_H Y })_{d-1} =0$.
\par
By Theorem \ref{HH} we know that
\begin{equation} \label{numeropuntisuspazio}
\dim (I_{ L_{x+1}+ \dots + L_{ e}})_{d-1} = {d+3 \choose 4} -  d(e-x) =x .
\end{equation}

Moreover, since the scheme $ \Pi+L_{x+1}+ \dots + L_{ e}$ has
$e-x $ lines,
and it is easy to show that
$$e-x={{d(d+2)(d+4)}\over 24} \geq
 \left \lceil {{{(d-1)+4 \choose 4} -  {(d-1)+2 \choose 2} }\over {(d-1)+1} } \right \rceil =  \left \lceil{{(d-1) (d+1)(d+6)} \over 24}  \right \rceil,
$$
then, by  the inductive hypothesis, we get
\begin{equation} \label{nullasuspazio}
\dim (I_{ \Pi+L_{x+1}+ \dots + L_{ e}})_{d-1}=0.
\end{equation}

Now we apply Lemma \ref{AggiungerePuntiSuSpazioLineare}; by (\ref{numeropuntisuspazio}) and (\ref {nullasuspazio}) we have
\begin{equation} \label{residuo}
\dim (I_{Res_H Y })_{d-1}=0.
\end{equation}
Finally, by (\ref{traccia}), (\ref{residuo}) and  Lemma
\ref{Castelnuovo} (with $\delta =1$)  we get $\dim (I_{ Y})_d =
0$, and that completes the proof of our main theorem
  for $\mathbb P^4$.
 \end{proof}

\section{Some technical lemmata}\label{lemmas}

Although the base case for an inductive approach to our main theorem
  was relatively straightforward, this is not the
case for the inductive step.

One aspect is relatively clear.  We first specialize some lines and degenerate other pairs of lines and divide our calculation, via Castelnuovo, into a Residual scheme (which we can handle easily) and a Trace scheme in a lower dimensional projective space.  It is here that the difficulties take place.  The Trace scheme will consist of degenerate conics, points and lines.  Unfortunately, it is not always the case that generic collections of degenerate conics behave well with respect to postulational questions.  The following example makes that clear.

\begin{rem} If $C$ is a degenerate conic in $\mathbb P^3$ then imposing the passage though $C$ imposes 7 conditions on the cubics of $\mathbb P^3$.  One might then suspect that if $X$ is the union of three generic degenerate conics in $\mathbb P^3$ then $X$ would impose $3\cdot 7 = 21$ conditions on cubics.  I.e. there would not be a cubic surface through $X$, although there obviously is one.
\end{rem}

It is the existence of such examples that complicates the induction step.  In fact, to get around this difficulty, we have to consider (at the same time) several auxiliary families combining both specializations and degenerations of a scheme consisting of a collection of generic lines and points.

Note that the first two lemmata deal with such families of auxiliary schemes in $\mathbb P^3$.  These are needed to deal with the Trace scheme in $\mathbb P^4$ which occurs in the first inductive step from $\mathbb P^4$ to $\mathbb P^5$.   These two lemmata also serve to point out the kinds of families we will need for the remainder of the proof.

\begin{lem} \label{RetteIncrociateInP3a}  Let $d=2(4h+r+1)$, $h \in \Bbb N$, $r=0;1;3$, (that is, $d \equiv 0; 2; 4,$ mod 8). Let
$$ c=   \left\lfloor{ {d+3 \choose 4} \over d}  \right\rfloor ,
$$
and set
$$a = {d+3 \choose 4}-d  c ; \ \ \ \ \  \ \ b={ {{d+3 \choose 3}- a(2d+1)-c } \over {d+1} } .$$
Then
\begin{itemize}
\item[(i)]
$b$ is an integer; \par
\item[(ii)]
if $x = {d+1 \choose 3}-(a+b)(d-1)$ we have $0\leq x<c ;$ \par
\item[(iii)]
if $W \subset \Bbb P^3$ is the following scheme
$$W = C_1+ \dots + C_a +M_1+ \dots +M_b + P_1+ \dots +P_c
$$
(where the $C_i$  are generic degenerate conics, the  $M_i$ are generic lines, and the $P_i$ are generic points)
then $W$ gives the expected number of conditions to the forms of degree $d$, that is
$$ \dim (I_W)_d  = {d+3 \choose 3} - a(2d+1)-b(d+1)-c =0.
$$
\end{itemize}
\end{lem}

\begin{proof}
{\rm (i)} An easy computation, yields
\begin{itemize}\item
 for  $d=8h+2$ (that is for $r=0$),
$$c={1\over4} {{d+3}\choose 3} -{1\over2};
\ \  a={d \over 2} = 4h+1; \ \  \hbox {and so} \
b= 8h^2+h+1 ;\ \ \
$$
\item
 for  $d=8h+4$ (that is for $r=1$),
$$c={1\over4} {{d+3}\choose 3} -{3\over4}; \ \
a={3d \over 4} = 6h+3; \ \  \hbox {and so} \
b= 8h^2+h;
$$
\item
 for  $d=8h+8$ (that is for $r=3$),
$$c={1\over4} {{d+3}\choose 3} -{1\over4};  \ \
a={d \over 4} = 2h+1; \ \  \hbox {and so} \
b= 8h^2+17h+10.$$
\end{itemize}

{\rm (ii)} Using  (i) and direct computation, (ii) easily follows.

{\rm (iii)} Observe that
$$ {d+3 \choose 3} - a(2d+1)-b(d+1)-c
$$
$$={d+3 \choose 3} - a(2d+1)-{d+3 \choose 3}+ a(2d+1)+c-c=0.$$
Thus we have to prove that $ \dim (I_W)_d  =0.$\par
If $d=2$, that is, for $h=r=0$, we have $a= 1$, $b=1$, $c= 2$, and it is  easy to see that there are not quadrics containing the scheme $C_1+M_1+P_1+P_2$.
\par
Let $d>2$.
Let $L_{i,1},L_{i,2}$ be  the two lines which form the  degenerate conic  $C_i$, and let $Q$ be a smooth quadric surface.
Let $x$ be as in {\rm (ii)}  and let $\widetilde W$ be the scheme obtained from $W$ by  specializing $(c-x)$ of the $c$ simple points $P_i$  to generic points on Q and by  specializing the conics $C_i$ in such a way that the
 lines $L_{1,1},\dots ,L_{a,1}$ become  lines of the same ruling on $Q$, (the lines $L_{1,2},\dots ,L_{a,2}$ remain generic lines,  not lying on $Q$).

 $L_{i,2}$ meets $Q$ in the two points which are $(L_{i,1} \cap L_{i,2})$ and another, which we denote by $R_{i,2}$. In the same way, $M_i$ meets $Q$ in the two points $S_{i,1}$, $S_{i,2}$.
 We have
 $$
 Res_Q  {\widetilde W} =L_{1,2}+\dots + L_{a,2} +M_1+ \dots+M_b+   P_1+\dots+P_x          \subset \Bbb P^4 ,
$$
where the $L_{i,2}$ and the $ M_i$ are generic lines.
By Theorem \ref{HH} and the description of $x$ we get
$$\dim ( I_{Res_Q  {\widetilde W}})_{d-2} =  {d+1 \choose 3}- (a+b)(d-1)-x=0.
$$
Now consider $
Tr_Q  {\widetilde W} $, which is
$$
L_{1,1}+ \dots  +L_{a,1}+R_{1,2}+\dots + R_{a,2}+
S_{1,1}+S_{1,2}+ \dots+S_{b,1}+ S_{b,2}+ P_{x+1}+\dots+P_c .
 $$
Note that the points $R_{i,2}, (1 \leq i \leq a); S_{i,1},S_{i,2}, (1 \leq i \leq b); P_i, (x+1 \leq i \leq c)$ are generic points on $Q$ and the lines all come from the same ruling on $Q$, hence
$$\dim ( I_{Tr_Q  {\widetilde W}})_{d} = (d-a+1)(d+1) - a -2b-(c-x).$$
By a  direct computation,  we get
$\dim ( I_{Tr_Q  {\widetilde W}})_{d} =0$.

So by Lemma \ref {Castelnuovo}, with $n=3$ and $\delta=2$, the conclusion follows.

\end{proof}

\begin{lem} \label{RetteIncrociateInP3b}  Let $d \geq3$ be odd, or $d=8h+6$, $h \in \Bbb N$, (that is, $ d\equiv 1;3;5;6;7$,  mod 8). Let
$$c=  { {d+3 \choose 4}\over {d}}$$
and set
$$ b=  \left\lfloor{ {d+4 \choose 4} \over {d+1} } \right\rfloor - c- 2;  \ \ \ \ \
b^*=  \left\lceil{ {d+4 \choose 4} \over {d+1} } \right\rceil - c- 2.$$
Then
\begin{itemize}
\item[(i)]
$b >0$ and $c$ is an integer; \par
\item[(ii)]
if $x = {d+1 \choose 3}-b(d-1)$, then $0\leq x<c;$ \par

\item[(iii)]
 if  $W, W^*  \subset \Bbb P^3$ are the following schemes
$$W = C+ 2P+M_1+ \dots +M_b + P_1+ \dots +P_c,
$$
$$W^* = C+ 2P+M_1+ \dots +M_{b^*} + P_1+ \dots +P_c,
$$
(where $C= L_1+L_2$  is a degenerate conic, formed by the two lines $L_1$, $L_2$; where $2P$ is a double point with support in $P= L_1 \cap L_2$; where   the  $M_i$ are generic lines and the $P_i$ are generic points)
then $W$ and $W^*$ give the expected number of conditions to the forms of degree $d$, that is
$$ \dim (I_W)_d  = {d+3 \choose 3} - (2d+2)-b(d+1)-c ,
$$
$$ \hbox {and} \ \ \  \dim (I_{W^*})_d  =0 .
$$
\end{itemize}
\end{lem}

\begin{proof}
Computing directly it is easy to verify {\rm (i)} and {\rm (ii)}.

{\rm (iii)} Since the scheme $C+2P$ is a degeneration of two skew lines it imposes $2d+2$ conditions to forms of degree $d$ (see Lemma \ref{sundial}). It follows that
 $$ \dim (I_W)_d  \geq {d+3 \choose 3} - (2d+2)-b(d+1)-c .$$

Hence, it suffices to prove that $ \dim (I_W)_d  \leq {d+3 \choose 3} - (2d+2)-b(d+1)-c .$

Let $Q$ be a smooth quadric surface.
Let $x$ be defined as in {\rm (ii)} and let  $\widetilde W$ be the scheme obtained from $W$ by  specializing $(c-x)$ of the $c$ simple points $P_i$  onto Q and by  specializing the line $M_1$ and the conic $C$ in such a way that the
 lines $M_1$ and  $L_{1}$ become   lines of the same ruling on $Q$ (the line $L_{2}$ remain a generic line,  not lying on $Q$, while the point $P$ becomes a point lying on $Q$).
We have $L_{2} \cap Q = P+ R$, and set $M_i \cap Q = S_{i,1}+S_{i,2}$, ($2 \leq i \leq b$).
Then
 $$
 Res_Q  {\widetilde W} =L_{2}+M_2+ \dots+M_b+   P_1+\dots+P_x          \subset \Bbb P^4 ;
$$
$$
Tr_Q  {\widetilde W}
= L_{1}+ M_1+2P|_Q +R + S_{2,1}+S_{2,2}+ \dots+S_{b,1}+ S_{b,2}+ P_{x+1}+\dots+P_c .
 $$
 By Theorem \ref{HH}  we immediately get
$$\dim ( I_{Res_Q  {\widetilde W}})_{d-2} =  {d+1 \choose 3}- b(d-1)-x=0.
$$
Thinking of $Q$ as $\mathbb P^1 \times \mathbb P^1$, we see that the forms  of degree $d$ in the ideal of $L_{1}+M_1+2P|_Q$
are  curves of type $(d-2,d)$ in $\mathbb P^1 \times \mathbb P^1$ passing through $P$,
since $P$ already belongs to $L_1$.  With that observation,  it is easy to check that
$$\dim ( I_{Tr_Q  {\widetilde W}})_{d} =( d-1)(d+1)- 2 - 2(b-1)-c+x$$
$$=  {d+3 \choose 3} - (2d+2)-b(d+1)-c.    $$
So by Lemma \ref {Castelnuovo},  with $n=3$ and $\delta=2$,  it follows that
$$ \dim (I_W)_d  = {d+3 \choose 3} - (2d+2)-b(d+1)-c ,$$
and we are finished with the schemes $W$.
\par
\medskip
We now consider the schemes $W^*$. If  $b =b^*$ (i.e., if $d \equiv 5, 6$, mod 8), we have $W^*=W$.  In this case it is easy to verify that  the number
$${d+3 \choose 3} - (2d+2)-b(d+1)-c $$ is  zero and so we are done.

So we are left with the case
$b^*=b+1$.
Let  $\widetilde W^*$ be the scheme obtained from $W^*$ by  specializing   $(c-x)$ of the $c$ simple points $P_i$,  the lines $M_1$ and $M_2$ and the conic $C$ in such a way that the
 lines $M_1, M_2, L_1$ are   lines of the same ruling on $Q$, and the line $L_{2}$ remains a generic line  not lying on $Q$. Note that the point $P$ becomes a point of $Q$.

 Set $L_{2} \cap Q = P+ R$, and set $M_i \cap Q = S_{i,1}+S_{i,2}$, ($3 \leq i \leq b^*$).
 We have
 $$
 Res_Q  {\widetilde W} =L_{2}+M_3+ \dots+M_{b^*}+   P_1+\dots+P_x          \subset \Bbb P^3 .
$$
and
$$
Tr_Q  {\widetilde W}
 = L_{1}+ M_1+M_2+ 2P|_Q +R $$
 $$+
S_{3,1}+S_{3,2}+ \dots+S_{b^*,1}+ S_{b^*,2}+ P_{x+1}+\dots+P_c .
$$
 By Theorem \ref{HH}  we immediately get
$$\dim ( I_{Res_Q  {\widetilde W}})_{d-2} =  {d+1 \choose 3}- b(d-1)-x=0.
$$
Using the same reasoning as above,  it is easy to check that
$$\dim ( I_{Tr_Q  {\widetilde W}})_{d} =\max \left\{0; ( d-2)(d+1)- 2 - 2(b-1)-c+x\right\}=0.    $$
By Lemma \ref {Castelnuovo}, with $n=3$ and $\delta=2$,  it follows that $ \dim (I_{W^*})_d  = 0 .$
\par
\end{proof}

We now formalize what we did in these last lemmata.
\medskip

 Let $n,d, a, b, c, \in \Bbb N$,  $n \geq 3$, $d >0, $ $a+b \leq d-1$,  and let
 $$t=  \left\lfloor{ {d+n \choose n} \over {d+1} } \right\rfloor; \ \ \ \ \ \ \ t^*=  \left\lceil{ {d+n \choose n} \over {d+1} } \right\rceil.$$
 Let $c \leq t-2(a+b)$,  $c^* \geq t^*-2(a+b)$.
 Let $ \widehat C_i$ be a 3-dimensional sundial (see Definition \ref{defsundial}), that is a generic degenerate conic with an embedded point, and let $M_i$ be a  generic line.

  Note that  $t \geq 2(d-1)$.
  \par
   \medskip
  Consider the following statements:
 \medskip
\begin{itemize}

  \item $S(n,d)$: {\it The scheme  \par \noindent
  $W(n,d) = \widehat C_1 + \dots +\widehat C_{d-1}  + M_1+\dots+M_{t-2(d-1)} \subset \Bbb P^n,$
   \par \noindent
   imposes the expected number of conditions to  forms of degree $d$, that is:
     \par \noindent
 $ \dim (I_{W(n,d)})_d= {d+n \choose n} - (2d+2)(d-1) - (d+1)(t-2(d-1))
 $ \par
 $=  {d+n \choose n} - t (d+1);  $}

 \medskip

\item $S^*(n,d)$:  {\it The scheme  \par \noindent
  $W^*(n,d) = \widehat C_1 + \dots +\widehat C_{d-1}  + M_1+\dots+M_{t^*-2(d-1)} \subset \Bbb P^n,$
  \par \noindent
   imposes the expected number of conditions to  forms of degree $d$, that is:
     \par \noindent
    $ \dim (I_{W^*(n,d)})_d= 0.  $}

 \medskip

  \item  $S(n,d;a,b,c)$:  {\it The scheme  \par \noindent
  $W(n,d;a,b,c) = \widehat C_1 + \dots +\widehat C_{a}  +  D_1+\dots+D_{b}+R_1+ \dots +R_b+M_1+\dots+M_{c} \subset \Bbb P^n,$
 \par \noindent
where the $D_i$ are  generic degenerate conics, and the $R_i$
are generic points, imposes the expected number of conditions to  forms of degree $d$,
 that is:
  \par \noindent
 $ \dim (I_{W(n,d;a,b,c)})_d= {d+n \choose n} - (2a+2b+c) (d+1).
 $}

  \medskip
  \item  $S^*(n,d;a,b,c^*)$:  {\it The scheme  \par \noindent
  $W^*(n,d;a,b,c^*) = \widehat C_1 + \dots +\widehat C_{a}  +  D_1+\dots+D_{b}+R_1+ \dots +R_b+M_1+\dots+M_{c^*} \subset \Bbb P^n,$
 \par \noindent
where the $D_i$ are  generic degenerate conics, and the $R_i$
are generic points, imposes the expected number of conditions to  forms of degree $d$,
 that is:
  \par \noindent
 $ \dim (I_{W(n,d;a,b,c^*)})_d= 0.
 $}

\end{itemize}

 \begin{lem} \label{DegenerareRette} Notation as above,
 \begin{itemize}
\item[(i)] if $S(n,d)$ holds, then $S(n,d;a,b,c)$ holds;

\item[(ii)] if $S^*(n,d)$ holds, then $S^*(n,d;a,b,c^*)$ holds.

  \end{itemize}
 \end{lem}
 \begin{proof}
 A degenerate conic with an embedded point is either a degeneration of two generic lines, or a specialization of a scheme which is the union of a degenerate conic and a simple generic point.
 Then by the semicontinuity of the Hilbert function, and
since a line imposes at most $d+1$ conditions to the forms of degree $d$, we get (i).

(ii) immediately follows from the semicontinuity of the Hilbert function.

 \end{proof}

\begin{lem} \label{S(4,d)} Notation as above,  let
 $$t=  \left\lfloor{ {d+4 \choose 4} \over {d+1} } \right\rfloor; \ \ \ \ \ \ \ t^*=  \left\lceil{ {d+4 \choose 4} \over {d+1} } \right\rceil.$$
Then $S(4,d)$ and  $S^*(4,d)$ hold, that is
 $$ \dim (I_{W(4,d)})_d=   {d+4 \choose 4} - t (d+1)  \ \  \ \  \hbox{and} \ \ \ \dim (I_{W^*(4,d)})_d= 0,$$
where
  $$W(4,d) = \widehat C_1 + \dots +\widehat C_{d-1}  + M_1+\dots+M_{t-2(d-1)} \subset \Bbb P^4,$$
  $$W^*(4,d) = \widehat C_1 + \dots +\widehat C_{d-1}  + M_1+\dots+M_{t^*-2(d-1)} \subset \Bbb P^4,$$
and $ \widehat C_i = C_i+2P_i|_{H_i}=L_{i,1}+L_{i,2}+2P_i|_{H_i}$.
\end{lem}

\begin{proof} By induction on $d$. For $d=1$ both conclusions follows from Theorem \ref {HH}. \par
Let $d>1$.

We consider two cases: \par
{\it Case 1}: $d=2(4h+r+1)$, $h \in \Bbb N$, $r=0;1;3$, (that is, $d\equiv 0; 2; 4$, mod 8).
In this case $t=t^*$, and we will prove that $\dim (I_{W(4,d)})_d=0$.
Consider
 $$c=  \left\lfloor{ {d+3 \choose 4} \over {d} } \right\rfloor   \ \ \ \hbox {and}  \ \ \
a = {d+3 \choose  4}  -  dc.$$
Note that:
 \begin{itemize}\item
 for  $d=8h+2$ (that is for $r=0$):
$$
c={1\over4} {{d+3}\choose 3} -{1\over2}; \ \ \ a={d \over 2};
$$
\item
 for  $d=8h+4$ (that is for $r=1$):
$$
c={1\over4} {{d+3}\choose 3} -{3\over4}; \ \ \ a={3d \over 4};
$$
\item
 for  $d=8h+8$ (that is for $r=3$):
$$
c={1\over4} {{d+3}\choose 3} -{1\over4}; \ \ \ a={d \over 4}.
$$
\end{itemize}
It is easy to check that
$$ 1 \leq a \leq d-1; \ \ \ \ \ \  0 \leq t-2a-c \leq t-2(d-1)
.$$

Let  $H \simeq \Bbb P^3$ be a generic hyperplane.
Let  $W_s(4,d)  $ be the scheme obtained from $W(4,d)$  by  specializing
 $t-2a-c$ lines   $M_1,\dots , M_{t-2a-c}$ into $H$
and by  specializing  $a$ degenerate conics  $\widehat C_1, \dots , \widehat C_{a}$,  in such a way that
$L_{i,1}+L_{i,2} \subset H$, but  $2P_i|_{H_i} \not\subset H $, for $1 \leq i \leq a$.

 So
$$
 Res_H  {W_s(4,d)} =P_1 + \dots +P_{a} + \widehat C_{a+1} + \dots +\widehat C_{d-1}  $$
 $$+ M_{t-2a-c+1}+\dots+M_{t-2(d-1)} \subset \Bbb P^4,
 $$
where $P_1, \dots ,P_{a}$ are generic points lying on $H$;
 $$
 Tr_H  {W_s(4,d)} =C_1 + \dots +C_{a} + R_{{a+1} ,1}+R_{{a+1} ,2} + \dots + R_{{d-1} ,1}+ R_{{d-1} ,2}
 +$$
 $$ +M_1+\dots + M_{t-2a-c} +S_{t-2a-c+1}+\dots+S_{t-2(d-1)} \subset \Bbb P^3,
 $$
where $R_{{i} ,1}+R_{{i} ,2} =  \widehat C_i \cap H = L_{i,1}\cap H+L_{i,2}\cap H$ and $S_{i}= M_i  \cap H$.

By Lemma \ref{RetteIncrociateInP3a},
since the $R_{{i} ,j}$  and the $S_i$ are
$$2(d-1-a)+ (t-2(d-1)-t+2a+c) = c $$ generic points, and $t-2a-c =b$ ($b$ as in  Lemma \ref{RetteIncrociateInP3a}), we get
$$\dim (I_{ Tr_H  {W_s(4,d)}}  )_d =0.$$
If we can prove that $\dim (I_{ Res_H  {W_s(4,d)}}  )_{d-1} =0$ then, by Lemma \ref{Castelnuovo},
with $\delta=1$, we are done.
 If $d=2$,
we have $a=1$, $c=2$  and
$$
 Res_H  {W_s(4,d)} =P_1   + M_{2}+\dots+M_{3} \subset \Bbb P^4.
 $$
Clearly  $\dim (I_{ Res_H  {W_s(4,d)}}  )_1 =0.$ \par
Now let $d>2$ and set
$$X= \widehat C_{a+1} + \dots +\widehat C_{d-1}  + M_{t-2a-c+1}+\dots+M_{t-2(d-1)} \subset \Bbb P^4,
 $$
 (hence $ Res_H  {W_s(4,d)} = X +P_1 + \dots +P_{a} $).
So
$X$ is the union of $d-1-a$ degenerate conics with an embedded point and $2a+c-2(d-1)$ lines.
The first step here is to show that $X$ imposes the right number of conditions to the forms of degre $d-1$.

By the induction hypothesis we have that $S(4,d-1)$ holds.  Since $d-1-a \leq d-3$,
and
$$X=
\widehat C_{a+1} + \dots +\widehat C_{d-1}  + M_{t-2a-c+1}+\dots+M_{t-2(d-1)}$$
is a
$$
W(4,d-1; d-1-a,0,2a+c-2(d-1)),
 $$
it follows from Lemma \ref{DegenerareRette} (i) that $X$ imposes independent conditions to the forms of degree $d-1$.
Thus
$$\dim (I_{X })_{d-1} =  {d-1+4 \choose 4}- d (2(d-1-a) + 2a+c-2(d-1) ) $$
$$=  {d+3 \choose 4}- dc =a.$$

To finish the argument we apply Lemma \ref{AggiungerePuntiSuSpazioLineare}. This requires us to prove that
$\dim  (I_{X +H})_{d-1}=0$. But
$$\dim  (I_{X +H})_{d-1}=\dim  (I_{X})_{d-2}.$$
For $d=2$, we obviously have
$\dim  (I_{X})_{d-2}=0.$

For $d>2$, by the inductive hypothesis $S^*(4,d-2)$ holds.
Since the parameters of $X$ (perhaps with fewer lines) satisfy the restrictions necessary to use Lemma \ref{DegenerareRette} (ii), we get that
$$S^*(4,d-2; d-1-a, 0, 2a+c-2(d-1))$$
 holds, that is, $\dim  (I_{X})_{d-2}=0.$

So, by Lemma \ref{AggiungerePuntiSuSpazioLineare}, we have
$$\dim (I_{ X+P_1 + \dots +P_{a} })_{d-1}  = \dim (I_{Res_H  {W_s(4,d)} })_{d-1} = 0,
$$
and we are done.
\par
\medskip

{\it Case 2}: $d$ odd, or $d=8h+6$, $h \in \Bbb N$, (that is, $ d= 1;3;5;6;7, $ mod  8).
Let
$$c=  { {d+3 \choose 4}\over {d}};  \ \ \ \ \  b= t - c- 2;\ \ \ \ \    b^*= t^* - c- 2,$$
(note that c is an integer).
It is easy to check that
$$0 <  b \leq t-2(d-1) \ \ \hbox{ and } \ \   0 <  b^* \leq t^*-2(d-1).$$

Let  $W_s(4,d)  $ be the scheme obtained from $W(4,d)$  by  specializing
the $b$ lines   $M_1,\dots , M_{b}$ and $\widehat C_{d-1}$
into a hyperplane $H \simeq \Bbb P^3$.

Let  $W^*_s(4,d)  $ be the scheme obtained from $W^*(4,d)$  by  specializing
into  $H$ the lines $M_1,\dots , M_{b^*}$ and the degenerate conic with an embedded point $\widehat C_{d-1}$.
We have
$$ Res_H  {W_s(4,d)} = \widehat C_1 + \dots +\widehat C_{d-2} + M_{b+1}+\dots+M_{t-2(d-1)} \subset \Bbb P^4,
 $$
 $$ Res_H  {W^*_s(4,d)} = \widehat C_1 + \dots +\widehat C_{d-2} + M_{b^*+1}+\dots+M_{t^*-2(d-1)} \subset \Bbb P^4,
 $$
 that is, both $Res_H  {W_s(4,d)} $ and $Res_H  {W^*_s(4,d)}$ are the union of $d-2$ degenerate conics with an embedded point and $c-2d+4$ lines.
 By the inductive hypothesis we immediately get
$$\dim (I_{Res_H  {W_s(4,d)}  })_{d-1} =$$

$$=\dim (I_{Res_H  {W^*_s(4,d)} })_{d-1} = {d+3 \choose 4} - d(
2(d-2)+c-2d+4)=0.
$$
Now we consider the traces:
$$
 Tr_H  {W_s(4,d)} =R_{{1} ,1}+R_{{1} ,2} + \dots + R_{{d-2} ,1}+ R_{{d-2} ,2}   + \widehat C_{d-2} +
 $$
 $$ +M_1+\dots + M_{b} +S_{b+1}+\dots+S_{t-2(d-1)} \subset \Bbb P^3,
 $$
 $$
 Tr_H  {W^*_s(4,d)} =R_{{1} ,1}+R_{{1} ,2} + \dots + R_{{d-2} ,1}+ R_{{d-2} ,2}   + \widehat C_{d-2} +
 $$
 $$ +M_1+\dots + M_{b^*} +S_{b^*+1}+\dots+S_{t^*-2(d-1)} \subset \Bbb P^3,
 $$
where $R_{{i} ,1}+R_{{i} ,2} =  \widehat C_i \cap H$,  and $S_{i}= M_i  \cap H$.

$Tr_H  {W_s(4,d)}$ is the union of $2(d-2)+c+4-2d= c$ simple generic  points, a  degenerate conic with an embedded point, and $b$ lines. So, by Lemma
\ref{RetteIncrociateInP3b} we get
$$ \dim (I_{ Tr_H  {W_s(4,d)} })_d  = {d+3 \choose 3} - (2d+2)-b(d+1)-c .
$$
Thus, by Lemma \ref{Castelnuovo}, with  $\delta=1$, we have
$$ \dim (I_  {W_s(4,d)} )_d  \leq  {d+3 \choose 3} - (2d+2)-b(d+1)-c = {d+4 \choose 4} - t(d+1).
$$
Since $  \dim (I_  {W(4,d)} )_d \leq \dim (I_  {W_s(4,d)} )_d$ and ${d+4 \choose 4} - t(d+1)$ is the expected dimension for $(I_  {W(4,d)} )_d$, we have
$  \dim (I_  {W(4,d)} )_d = {d+4 \choose 4} - t(d+1)$.
\par
Finally,
$Tr_H  {W^*_s(4,d)}$ is the union of $ c$ simple generic  points, one  degenerate conic with an embedded point, and $b^*$ lines.
So, by Lemma  \ref{RetteIncrociateInP3b} we get
$$ \dim (I_{ Tr_H  {W^*_s(4,d)} })_d  =0,
$$
and by Lemma \ref{Castelnuovo}, with  $\delta=1$, the conclusion follows.

\end{proof}

\begin{lem} \label{S(n,d)} Let $n \geq4$, $d\geq1$,
 $$t=  \left\lfloor{ {d+n \choose n} \over {d+1} } \right\rfloor; \ \ \ \ \ \ \ t^*=  \left\lceil{ {d+n \choose n} \over {d+1} } \right\rceil.$$
Then $S(n,d)$ and  $S^*(n,d)$ hold, that is
 $$ \dim (I_{W(n,d)})_d=   {d+n \choose n} - t (d+1);  \ \ \ \ \dim (I_{W^*(n,d)})_d= 0,$$
where
  $$W(n,d) = \widehat C_1 + \dots +\widehat C_{d-1}  + M_1+\dots+M_{t-2(d-1)} \subset \Bbb P^n,$$
  $$W^*(n,d) = \widehat C_1 + \dots +\widehat C_{d-1}  + M_1+\dots+M_{t^*-2(d-1)} \subset \Bbb P^n,$$
and $ \widehat C_i = C_i+2P_i|_{H_i}=L_{i,1}+L_{i,2}+2P_i|_{H_i}$.
\end{lem}

\begin{proof} By induction on $n+d$. The case $d=1$ follows from Theorem \ref {HH}. For $n=4$, see Lemma \ref{S(4,d)}.
\par
Let $n+d>6$, $d>1$, $n>4$. Let
$$ a=  {d+n-1 \choose n} -  d \left\lfloor{ {d+n-1 \choose n} \over {d} } \right\rfloor \ \ \ \hbox{  and } \ \ \  c=\left\lfloor{ {d+n-1 \choose n} \over {d} } \right\rfloor
.$$
Note that, by a direct computation, we have $$0 \leq a \leq d-1 \ \ \ \hbox{ and } \ \ \ \  a \leq c  \leq t-2(d-1).$$
\par
Let  $W_s(n,d)  $ be the scheme obtained from $W(n,d)$  by  specializing,
into a generic hyperplane $H \simeq \Bbb P^{n-1}$, the $d-1-a$
degenerate conics with an embedded point $\widehat C_{a+1}$,  $\dots, \widehat C_{d-1}$
 and the $ t-2(d-1)- c $ lines $M_{c+1}, \dots ,M_{t-2(d-1)}$. We further  specialize
 the  $a$ degenerate conics  $\widehat C_1, \dots , \widehat C_{a}$, in such a way that $L_{i,1}+L_{i,2}\subset H$, but
 $2P_i|_{H_i} \not\subset H $, for $1 \leq i \leq a$.

 \par
 Analogously, let  $W^*_s(n,d)  $ be the scheme obtained from $W^*(n,d)$  by  specializing,
 into a generic hyperplane $H \simeq \Bbb P^{n-1}$,
the degenerate conics with an embedded point $\widehat C_{a+1}, \dots ,\widehat C_{d-1}$,
 and  the $ t^*-2(d-1)- c $ lines $M_{c+1}, \dots , M_{t^*-2(d-1)}$.
We further  specialize
 the  $a$ degenerate conics  $\widehat C_1, \dots , \widehat  C_{a}$, in such a way that $L_{i,1}+L_{i,2}\subset H$, but
 $2P_i|_{H_i} \not\subset H $.

From these specializations we have
$$
 Res_H  {W_s(n,d)} = Res_H  {W^*_s(n,d)} =P_1 + \dots +P_{a}  + M_{1}+\dots+M_{c} \subset \Bbb P^n,
 $$
where $P_1, \dots ,P_{a}$ are generic points of $H$;
 $$
 Tr_H  {W_s(n,d)} =
 $$
$$
 C_1 + \dots +C_{a} +  \widehat C_{a+1} + \dots +\widehat C_{d-1}  +S_{1}+\dots+S_{c} +M_{c+1}+\dots + M_{t-2(d-1)} \subset \Bbb P^{n-1},
 $$
and
  $$
 Tr_H  {W^*_s(n,d)} =
 $$
$$
 C_1 + \dots +C_{a} +  \widehat C_{a+1} + \dots +\widehat C_{d-1}  +S_{1}+\dots+S_{c} +M_{c+1}+\dots + M_{t^*-2(d-1)} \subset \Bbb P^{n-1},
 $$
where  $S_{i}= M_i  \cap H$.

Consider the schemes
 $$
X=  Tr_H  {W_s(n,d)} - (S_{a+1}+\dots+S_{c})$$
$$=  \widehat C_{a+1} + \dots +\widehat C_{d-1} +C_1 + \dots +C_{a}  +S_{1}+\dots+S_{a} +M_{c+1}+\dots + M_{t-2(d-1)} \subset \Bbb P^{n-1},
 $$
and
  $$
X^*=  Tr_H  {W^*_s(n,d)} - (S_{a+1}+\dots+S_{c})$$
$$=  \widehat C_{a+1} + \dots +\widehat C_{d-1} +C_1 + \dots +C_{a}  +S_{1}+\dots+S_{a} +M_{c+1}+\dots + M_{t^*-2(d-1)} \subset \Bbb P^{n-1}.
 $$
By the inductive hypothesis, $S(n-1,d)$ holds.  By a direct
computation we check that
$$t-c \leq t^* -c \leq   \left\lfloor{ {d+n-1 \choose d} \over {d+1} } \right\rfloor .$$

Hence,   by Lemma \ref{DegenerareRette},
we have that  $S(n-1,d; d-1-a, a, t-2(d-1)-c)$ and $S(n-1,d; d-1-a, a, t^*-2(d-1)-c)$ hold.

It follows that
$\dim (I_X)_{d}$ and $\dim (I_{X^*})_{d}$ are as expected, that is,
$$\dim (I_X)_{d} = {d+n-1 \choose {n-1}}  - (d+1)(2(d-1)+t-2(d-1)-c)
$$
$$=  {d+n-1 \choose {n-1}}  - (d+1)(t-c)=  {d+n \choose {n}} -t(d+1)+c-a,
$$
and
$$\dim (I_{X^*})_{d} = {d+n-1 \choose {n-1}}  - (d+1)(2(d-1)+t^*-2(d-1)-c)
$$
$$=  {d+n-1 \choose {n-1}}  - (d+1)(t^*-c)=  {d+n \choose {n}} -t^*(d+1)+c-a.
$$
Now, since $S_{a+1}, \dots,S_{c}$ are generic points and  $ {d+n \choose {n}} -t^*(d+1) \leq 0 $, it follows that
$$\dim (I_{ Tr_H  {W_s(n,d)} })_{d} ={d+n \choose {n}} -t(d+1),
$$
and
$$\dim (I_{ Tr_H  {W^*_s(n,d)} })_{d} =\max \left\{ 0 ; {d+n \choose {n}} -t^*(d+1) \right\} =0.
$$
If we prove that $\dim (I_{ Res_H  {W_s(n,d)} })_{d-1} = \dim (I_{ Res_H  {W^*_s(n,d)} })_{d-1} =0$ the, by Lemma \ref{Castelnuovo} with  $\delta=1$, we are done.

Recall that
$$
 Res_H  {W_s(n,d)} = Res_H  {W^*_s(n,d)} =P_1 + \dots +P_{a}  + M_{1}+\dots+M_{c} \subset \Bbb P^n,
 $$
where $P_1, \dots ,P_{a} $ are generic points in $H$. By Lemma \ref{AggiungerePuntiSuSpazioLineare} it sufficies to prove that
$\dim (I_{M_{1}+\dots+M_{c}  })_{d-1} =a$ and $\dim (I_{M_{1}+\dots+M_{c}+H})_{d-1} =0$.

By Theorem \ref{HH} we immediately get
$$\dim (I_{M_{1}+\dots+M_{c}  })_{d-1} ={d+n-1 \choose {n}}  - dc= a.
$$
Moreover, since  $\dim (I_{M_{1}+\dots+M_{c}+H})_{d-1} =\dim (I_{M_{1}+\dots+M_{c}})_{d-2}$,
by Theorem \ref{HH} we have
$$\dim (I_{M_{1}+\dots+M_{c}})_{d-2} = \max \left\{ 0; {d+n-2 \choose
{n}}  - (d-1) c \right\}=0,$$ and the conclusion follows.

\end{proof}

\section{The general case}\label{generalsection}

Having collected all the preliminary lemmata necessary, we are
ready to prove the main theorem of the paper.

\begin{thm} \label{TeoremaInPn} Let $n,d \in \Bbb N$, $n  \geq 4$, $d \geq1$. Let $\Pi \subset \Bbb P^n$ be a plane, and let $L_1, \dots , L_s \subset \Bbb P^n$ be $s$ generic lines.
If
 $$X= \Pi +  L_1+ \dots + L_s \subset \Bbb P^n$$
then
 $$
 \dim (I_{X})_d = \max \left \{ {d+n \choose n} -  {d+2 \choose 2}  -s(d+1), 0 \right
 \},
 $$
or equivalently $X$ has bipolynomial Hilbert function.
 \end{thm}

 \begin{proof}
 We proceed by induction on $n+d$. The result is obvious for $d=1$ and any $n$, while for $n=4$ see Theorem  \ref{TeoremaInP4}.  \par
 Let $d>1$, $n>4$.
 By  Lemma \ref{BastaProvarePers=e,e*} it suffices to prove the theorem for $s=e$ and $s=e^*$, where
 $$e= \left \lfloor   {{{d+n \choose n} -  {d+2 \choose 2} }\over {d+1} }\right \rfloor ; \ \ \ \
 e^*= \left \lceil   {{{d+n \choose n} -  {d+2 \choose 2} }\over {d+1} }\right \rceil .
 $$
Let
$$e_\rho = \left \lfloor   {{{d+n-1 \choose n} -  {d+1 \choose 2} }\over {d} }\right \rfloor  ;   \ \ \ \
\rho=  {{d+n-1 \choose n} -  {d+1 \choose 2} }- e_\rho d ;$$   $$e_T = s - e_\rho - 2\rho, \ \ \ \ \  (s=e, e^*).
$$
It is a direct computation to check that
$e - e_\rho - 2\rho \geq 0$.

Let  $\widehat C_i $ be the degenerate conic with an embedded
point obtained by degenerating the  lines $L_{i},  L_{i+1}$, $1
\leq i \leq \rho$ as in Lemma \ref{sundial} with $m=1$. By abuse
of notation, we write $\widehat C_i $ as $L_{i}+  L_{i+1}+2P_i|
_{H_i}$, (recall that $H_i \simeq \Bbb P^3 $ is a generic linear
space through $P_i$). Let $H \simeq \Bbb P^{n-1}$ be a generic
hyperplane. Now specialize  $\widehat C_1, \dots, \widehat C_\rho
$ in such a way that  $L_{i}+  L_{i+1} \subset H$ and $2P_i|_{H_i}
\not\subset H $, and specialize the $e_T$ lines $L_{2\rho+1},
\dots L_{2\rho+e_T}$  into $H$ and denote by $Y $ the resulting
scheme. We have
$$
Res_H  { Y} = \Pi + P_1+ \dots  + P_\rho + L_{2\rho+e_T+1}+ \dots
+L_{s} \subset \Bbb P^{n}
$$
($P_1, \dots ,P_{\rho}$ are generic points of $H$),
$$
Tr_H  { Y}  =L + C_1+ \dots  + C_\rho + L_{2\rho+1}+ \dots +
L_{2\rho+e_T}+ P_{2\rho+e_T+1}+ \dots +P_{s} \subset \Bbb P^{n-1}
$$
where $L =\Pi  \cap H$ and  $P_i = L_i \cap H$ , for $2\rho+e_T+1 \leq i \leq s$.

$Res_H  { Y} $ is the union of  one plane, $e_\rho$ lines and
$\rho$ generic points of $H$. By the inductive hypothesis we have
$$\dim (I_{ \Pi  + L_{2\rho+e_T+1}+ \dots +L_{s} })_{d-1} =\rho.
$$
Moreover
$$\dim (I_{ H+ \Pi  + L_{2\rho+e_T+1}+ \dots +L_{s} })_{d-1} =\dim (I_{  \Pi  + L_{2\rho+e_T+1}+ \dots +L_{s} })_{d-2}=0,
$$
(obvious, for $d=1,2$; by induction, for $d>2$).
Hence by
Lemma \ref{AggiungerePuntiSuSpazioLineare} we get
$$\dim (I_{ Res_H  { Y}})_{d-1} =0.
$$

$Tr_H  { Y} $  is the union of  $\rho$ degenerate conics, $e_T +1$
lines, and $e_\rho$ generic points. We will compute $\dim (I_{
Tr_H  { Y}})_{d} $ by using Lemma \ref{DegenerareRette} and Lemma
\ref{S(n,d)}. We have to check that  $\rho \leq d-1$ and $e_\rho
\leq \rho$. The first inequality is obvious, and it is not
difficult  to verify the other one. So we get
$$\dim (I_{ Tr_H  { Y}})_{d} = \max \left\{  0; {d+n-1 \choose d}    - (d+1)(2\rho+e_T+1)-e_\rho\right\} ,
$$
and from here
$$\dim (I_{ Tr_H  { Y}})_{d} = {d+n \choose n} -  {d+2 \choose 2}  -s(d+1), \ \ \ \  for \ \ s=e ;$$
$$\dim (I_{ Tr_H  { Y}})_{d} =0 \ \ \ \  for \ \ s=e^*  .$$
 The conclusion now follows from Theorem \ref{Castelnuovo}, with  $\delta=1$.

\end{proof}

\section{Applications}\label{applicationsection}

We now mention two applications of Theorem \ref{TeoremaInPn}. The
first is to a very classical problem concerning the existence of
rational normal curves having prescribed intersections with
various dimensional linear subspaces of $\PP^n$.  For example, the
classical Theorem of Castelnuovo which asserts that there exists a
unique rational normal curve through $n+3$ generic points of
$\PP^n$, is the kind of result we have in mind.

The second application is to writing polynomials in several
variables in a simple form. For example, the classical theorem which says that
in $S=\mathbb C[x_0, \dots, x_n]$ every quadratic form is a sum of
at most $n+1$ squares of linear forms, is the kind of  theorem we
intend.

\medskip
\noindent {\bf Rational normal curves.} The problem of deciding
whether or not there exists a rational normal curve with
prescribed intersections with generic configurations of linear spaces, is
well known and, in general, unsolved. Various results and
applications of answers to this problem can be found in
\cite{CaCat07} and \cite{CaCat09}.

Of particular importance in such questions is the Hilbert function of the resulting configuration of linear spaces.
It is for this reason that the results of this paper can be applied to such a problem.

To illustrate the relationships we will look at the following
special problem (left open in \cite{CaCat09}): consider in
$\mathbb P^4$, $P_1,P_2,P_3$ generic points, $L_1,L_2$  generic
lines and  $\pi$ a generic plane. Does there exist a rational
normal curve $\mathcal{C}$ in $\mathbb  P^4$ such that:

(i) $\mathcal{C}$ passes through the $P_i$ ($i=1,2,3$);

(ii) $\deg(\mathcal{C}\cap L_i)\geq 2$ for $i=1,2$;

(iii) $\deg(\mathcal{C}\cap\pi)\geq 3$.

 An expected answer is described in \cite{CaCat09} and can be obtained by arguing as follows:
inside  the 21 dimensional parameter space for rational normal
curves in $\PP^4$ it is  expected that those satisfying the
conditions enumerated above form a subvariety of codimension $20$. In other words,
we expect that there is a rational normal curve in $\mathbb P^4$
satisfying the conditions above.

To see that this is not the case  we consider the schemes
$${X}=P_1+P_2+L_1 +L_2+\pi,$$
$${Y}={X}+P_3.$$

Using Theorem \ref{TeoremaInPn} we know that $\dim (I_{ X})_2 =1$
and $\dim (I_{ Y})_2 =0$. If $\mathcal{C}$ existed, then
$Q\supset{X}$ would imply $Q\supset\mathcal{C}$ by a standard
Bezout type argument, and so we get $Q\supset {Y}$, a
contradiction.

\medskip
\noindent{\bf Polynomial decompositions.} We consider the rings
$S=\mathbb{C}[x_0,\ldots,x_n]$ and $T=\mathbb{C}[y_0,\ldots,y_n]$,
and we denote by $S_d$ and $T_d$ their homogeneous pieces of
degree $d$. We consider $T$ as an $S$-module by letting the action
of $x_i$ on $T$ be that of partial differentiation with respect to
$y_i$. We also use some basic notions about apolarity (for more on
this see \cite{Ge,IaKa}).

Let $I\subset S$ be a subset and denote by $I^\perp\subset T$ the
submodule of $T$ annihilated by every element of $I$. If $I$ is an
homogeneous ideal, we recall that $(I_d)^\perp=(I^\perp)_d$.

Given linear forms $a,b,c,l_i,m_i\in T_1,i=1,\ldots,s,$ one can ask
the following question $(\star)$:
\begin{quote}{\em
For which values of $d$ is it true that any form $f\in T_d$ can be
written as
$$f(y_0,\ldots,y_n)=f_1(l_1,m_1)+\ldots+f_s(l_s,m_s)+g(a,b,c)$$ for
suitable forms $f_i$ and $g$ of degree $d$?}
\end{quote}

\noindent More precisely, we ask whether the following vector space equality
holds:
$$T_d=\left(\mathbb{C}[l_1,m_1]\right)_d+\ldots
+\left(\mathbb{C}[l_s,m_s]\right)_d+\left(\mathbb{C}[a,b,c]\right)_d,$$
where $\left(\mathbb{C}[l_i,m_i]\right)_d$, respectively
$\left(\mathbb{C}[a,b,c]\right)_d$, is the degree $d$ part of the
subring of $T$ generated by the $l_i,m_i$'s for a fixed $i$,
respectively generated by $a,b$ and $c$. A more general question
can be considered as described in \cite{CarCatGer1}, but a
complete answer is not known. We now give a complete answer in the case of $(\star)$.

The connection with configurations of linear spaces is given by the
following results.

\begin{lem} Let $\Lambda\subset\PP^n$ be an $i$ dimensional
linear space having defining ideal $I$. Then, for any $d$, we have
the following:

$$ I_d^\perp=\left(\mathbb{C}[l_0,\ldots,l_i]\right)_d $$

where the linear forms $l_i\in T_1$ generate $I_1^\perp$.
\end{lem}
\begin{proof}
After a linear change of variables, we may assume
$$I=(x_0,\ldots,x_{n-i-1}).$$ As this is a monomial ideal the
conclusion follows by straightforward computations.
\end{proof}

\begin{prop}
Let $\Lambda=\Lambda_1+\ldots+\Lambda_s\subset\PP^n$ be a
configuration of linear spaces having defining ideal $I$ and such
that $\dim \Lambda_i=n_i$. Then, for any $d$, the following holds:
$$ I_d^\perp=\left(\mathbb{C}[l_{1,0},\ldots,l_{1,n_1}]\right)_d+\ldots +\left(\mathbb{C}[l_{s,0},\ldots,l_{s,n_s}]\right)_d$$
where the linear forms $l_{i,j}\in T_1$ are such that the degree
$1$ piece of $(l_{i,0},\ldots,l_{i,n_i})^\perp$ generates the
ideal of $\Lambda_i$.
\end{prop}
\begin{proof}
The proof follows readily from the previous lemma once we recall
that $(I\cap J)^\perp=I^\perp + J^\perp$.
\end{proof}

Now we can make clear the connection with question $(\star)$.
Given the linear forms $a,b,c,l_i,m_i\in T_1$ for $i=1,\ldots,s$,
we consider the ideal $I\subset S$ generated by the degree 1 piece
of $(a,b,c)^\perp$ and the ideals $I_i$ generated by the degree 1
pieces of $(l_i,m_i)^\perp, i=1,\ldots,s$. Note that $I\cap
I_1\cap\ldots\cap I_s$ is the ideal of the union of $s$ lines and
one plane in $\PP^n$. Denote this scheme by ${X}$. Now we can give
an answer to question $(\star)$ using Theorem \ref{TeoremaInPn}.

\begin{prop} With notation as above, we have: the values of $d$ answering question $(\star)$  are exactly the ones
for which $\dim (I_{ X})_d=0$.
\end{prop} \qed

\section{Final remarks}

Theorem \ref{TeoremaInPn} gives new evidence for the conjecture we
stated in the Introduction of the paper. As our conjecture deals
with generic configurations of linear spaces with non-intersecting
components, we would like to say something in case there are components which are forced to intersect.

Let $\Lambda=\bigcup\Lambda_i\subset\PP^n$ be a generic
configuration of linear spaces such that $m_i=\dim\Lambda_i\geq
m_j=\dim\Lambda_j$ if $i\geq j$. Then, there exist components of
$\Lambda$ which intersect if and only if $m_1+m_2\geq n$.

The first interesting case where generic configurations of linear
spaces have intersecting components occurs  in $\PP^3$ by taking lines
and at least one plane.

\begin{rem} Theorem \ref{TeoremaInPn} is not stated in $\PP^3$,
but it can easily be extended to include this case. If
$X=L_1+\ldots+L_s+\Pi\subset\PP^3$ we consider the exact sequence
\[0\rightarrow I_{L_1+\ldots +L_s}(-1)\rightarrow R \rightarrow R/I_X \rightarrow 0\]
where the first map is  multiplication by a linear form
defining $\Pi$. We can compute $HF(X,\cdot)$ by taking
dimensions in degree $d$ and obtain:
\[HF(X,d)={d+3\choose 3}-\max\left\{0,{d+2\choose 3}-sd\right\}\]
for $d>0$ and $HF(X,0)=1$. We also notice that
\[hp(X,d)={d+2\choose 2}+s(d+1)-s={d+3\choose 3}-{d+2\choose 3}+sd.\]
Thus $X$ has bipolynomial Hilbert function.
\end{rem}

Hence our conjecture holds for the union of generic lines and {\bf one}
plane even in $\PP^3$, where forced intersection appear. But, in
general, our conjecture is false for configurations of linear
spaces with intersecting components, as shown by the following
example.

\begin{ex}\label{notmaxrem}{
Consider $\Lambda\subset\PP^3$ a generic configuration of linear
spaces consisting of one line and three planes. By Derksen's
result in \cite{Derksen} we have $hp(\Lambda,1)=3$ but clearly no
plane containing $\Lambda$ exists. Hence,
\[HF(\Lambda,1)=4\neq\min\{hp(\PP^3,4)=4,hp(\Lambda,4)=3\}\] and
the Hilbert function is not bipolynomial.}
\end{ex}

We are not aware of any general result providing evidence for the
behavior of $HF(\Lambda,d)$ when the components of $\Lambda$ are
intersecting. We did, however, conduct experiments using the
computer algebra system CoCoA \cite{cocoa} and the results
obtained suggest the following:

\begin{quote}{\em let $\Lambda \subset\PP^n$ be a generic configuration of
linear spaces. There exists an integer $d(\Lambda)$ such that
\[HF(\Lambda,d)=hp(\PP^n,d), \mbox{ for } d\leq d(\Lambda)\]
and
\[HF(\Lambda,d)=hp(\Lambda,d), \mbox{ for } d> d(\Lambda).\]}
\end{quote}
This seems to be a reasonable possibility for the Hilbert function
of generic configurations of linear spaces (even with forced
intersections), but the evidence is still to sparse to call it a
conjecture.

\def\cprime{$'$}

\end{document}